\documentclass[12pt,leqno]{article}
\usepackage{amsmath, amsthm, amsfonts, amssymb, color}
\usepackage{mathrsfs}
\usepackage[notcite,notref]{showkeys}

\theoremstyle{definition}

\newcommand{\scr}[1]{\mathscr #1}
\definecolor{wco}{rgb}{0.5,0.2,0.3}

\numberwithin{equation}{section}

\newcommand{\ua}{\uparrow}

\title{{\bf Stochastic Heat Equations with Values in a Riemannian Manifold}\footnote{Supported in
 part by  NSFC (11401019, 11671035, 11371099) and DFG through  CRC 701}}
\author{
{\bf    Michael R\"{o}ckner$^{a)}$, Bo Wu$^{b,e)}$, Rongchan Zhu$^{a,c)}$, Xiangchan Zhu$^{a,d)}$}
\thanks{E-mail address: roeckner@math.uni-bielefeld.de(M.R\"{o}ckner), wubo@fudan.edu.cn(B.Wu),
zhurongchan@126.com(R.C.Zhu), zhuxiangchan@126.com(X.C.Zhu)}
\\
\footnotesize{  $^{a)}$ Department of Mathematics, University of Bielefeld, D-33615 Bielefeld, Germany}\\
\footnotesize {$^{b)}$ School  of Mathematical Sciences, Fudan
University, Shanghai 200433, China}\\
 \footnotesize{ $^{c)}$Department of Mathematics, Beijing Institute of Technology, Beijing 100081,  China}\\
\footnotesize{ $^{d)}$School of Science, Beijing Jiaotong University, Beijing 100044, China}
\\
 \footnotesize{ $^{e)}$ Institute for Applied Mathematics, University of Bonn,
Bonn 53115, Germany
}}

\date{}

\begin{document}
\maketitle
\def\R{\mathbb R} \def\EE{\mathbb E} \def\P{\mathbb P}\def\Z{\mathbb Z} \def\ff{\frac} \def\ss{\sqrt}
\def\H{\mathbb H}
\def\HH{\mathbf{H}}
\def\DD{\Delta} \def\vv{\varepsilon} \def\rr{\rho}
\def\<{\langle} \def\>{\rangle} \def\GG{\Gamma} \def\gg{\gamma}
\def\ll{\lambda} \def\LL{\Lambda} \def\nn{\nabla} \def\pp{\partial}
\def\dd{\text{\rm{d}}}
\def\Id{\text{\rm{Id}}}\def\loc{\text{\rm{loc}}} \def\bb{\beta} \def\aa{\alpha} \def\D{\scr D}
\def\E{\scr E} \def\si{\sigma} \def\ess{\text{\rm{ess}}}
\def\beg{\begin} \def\beq{\beg}  \def\F{\scr F}
\def\Ric{\text{\rm{Ric}}}
\def\Vol{\text{\rm{Vol}}}
\def\Var{\text{\rm{Var}}}
\def\Ent{\text{\rm{Ent}}}
\def\Hess{\text{\rm{Hess}}}\def\B{\scr B}
\def\e{\text{\rm{e}}} \def\ua{\underline a} \def\OO{\Omega} \def\b{\mathbf b}
\def\oo{\omega}     \def\tt{\tilde} \def\Ric{\text{\rm{Ric}}}
\def\cut{\text{\rm{cut}}} \def\P{\mathbb P} \def\K{\mathbb K}
\def\ifn{I_n(f^{\bigotimes n})}
\def\fff{f(x_1)\dots f(x_n)} \def\ifm{I_m(g^{\bigotimes m})} \def\ee{\varepsilon}
\def\C{\scr C}
\def\PP{\scr P}
\def\M{\scr M}\def\ll{\lambda}
\def\X{\scr X}
\def\T{\scr T}
\def\A{\mathbf A}
\def\LL{\scr L}\def\LLL{\Lambda}
\def\gap{\mathbf{gap}}
\def\div{\text{\rm div}}
\def\Lip{\text{\rm Lip}}
\def\dist{\text{\rm dist}}
\def\cut{\text{\rm cut}}
\def\supp{\text{\rm supp}}
\def\Cov{\text{\rm Cov}}
\def\Dom{\text{\rm Dom}}
\def\Cap{\text{\rm Cap}}\def\II{{\mathbb I}}\def\beq{\beg{equation}}
\def\sect{\text{\rm sect}}\def\H{\mathbb H}

\begin{abstract}The main result of this note is the existence of  martingale solutions to the stochastic heat equation (SHE) in a Riemannian manifold by using suitable Dirichlet forms
on the corresponding path/loop space. Moreover, we present some  characterizations of the lower bound of the Ricci curvature by functional inequalities of various associated  Dirichlet forms.
\end{abstract}


\section{Introduction}\label{sect1}
This work is motivated by Tadahisa Funaki's pioneering work \cite{Fun92} for regular noise and Martin Hairer's recent construction \cite{Hai16}  with singular noise of a
natural evolution on the loop space over a Riemannian manifold $(M,g)$.  Both consider the formal Langevin dynamics
associated to the energy
$$E(u)=\frac{1}{2}\int_{S^1} g_{u(x)}(\partial_xu(x),\partial_xu(x))\dd x,$$
for smooth functions $u:S^1\rightarrow M$. One would like to build a Markov process $u$ taking
values in loops over $M$ with
 invariant (even symmetrizing) measure formally given by $\exp(-2E(u))Du$. A natural way of interpreting $\exp(-2E(u))Du$ is to think of it as
the Brownian bridge measure on $M$. See \cite{AD99} for proofs that natural approximations of $\exp(-2E(u))Du$ do indeed converge to Wiener
measure on $C([0,1];M)$.

Processes with invariant (even symmetrizing) measure given by Wiener
measure on $C([0,1];M)$  were first constructed  in the nineties by using the Dirichlet form given by the Malliavin gradient on path and loop spaces over  Riemannian manifolds, see
\cite{DR92, ALR93}. In this case, we call the associated Dirichlet form \textbf{O-U Dirichlet form}. For an alternative approach, not based on Dirichlet forms, see \cite{Nor98}. After that there were several follow-up papers concentrating on non-compact Riemannian manifold,
see \cite{CW14, WW1}. In particular,  when $M=\mathbb{R}^d$ these processes  correspond to the Ornstein-Uhlenbeck
processes from Malliavin calculus. When $M=\mathbb{R}^d$  the stochastic heat equation also admits Wiener
measure as the invariant measure.
To construct the solution to the stochastic heat equation on Riemannian manifold, in \cite{Hai16} Martin Hairer  wrote the equation in  local coordinates informally as:
\begin{equation}\label{eq1.1}\dot{u}^\alpha=\partial_x^2u^\alpha+{\Gamma^\alpha}_{\beta\gamma}(u)\partial_xu^\beta\partial_xu^\gamma+\sigma_i^\alpha(u)\xi_i,\end{equation}
where Einstein’s convention of summation over repeated indices is applied and
${\Gamma^\alpha}_{\beta\gamma}$ are the Christoffel symbols for the Levi-Civita connection of $(M, g)$, $\sigma_i^\alpha$ are the local coordinates for the vector fields $\sigma_i$ on $M$
satisfying $g_u(h,\bar{h})=\sum_ig_u(h,\sigma_i)g_u(\bar{h},\sigma_i)$ for  $h,\bar{h}\in T_uM$, and $\xi_i$ is a collection of independent space-time white noises.  Equation \eqref{eq1.1} may be considered
as some kind of a multi-component version of the KPZ equation.
By regularity structure theory, recently developed in \cite{Hai14, BHZ16, CH16}, local well-posedness of \eqref{eq1.1} has been obtained in \cite{Hai16}.

In this note, we construct a new Dirichlet form ($L^2$\textbf{-Dirichlet form}) such that the associated Markov process solves the stochastic heat equation (SHE) with values in a Riemannian manifold.
Moreover, we obtain some new characterizations of the lower bound of the Ricci curvature in terms of ${\bf L^2}$-gradient and functional inequalities associated to the above Dirichlet form.
In addition, we also prove the logarithmic Sobolev inequality holds on the path space over a Riemannian manifold with lower bounded  Ricci curvature. As a consequence, for the process we have $L^2$-exponential
ergodicity, recurrent irreducibility and the strong law of large numbers.

In Sections 2 and 3 below, we present and discuss these results in detail and explain the framework. We also sketch some proofs.  The details of the proofs are contained in \cite{RWZZ17}.

\section{A Diffusion Process on Path Space}\label{sect2}

Throughout this article, suppose that $M$ is a complete and stochastically complete Riemannian manifold with dimension $d$, and $\rho$ be the Riemannian distance on $M$.
Fix $o\in M$ and $T>0$. The based path space $W_{o}(M)=\{\gamma\in C([0,1];M):\gamma(0)=o\}$, which is a Polish
space under the uniform distance
$$d_\infty(\gamma,\sigma):=\displaystyle\sup_{t\in[0,1]}\rho(\gamma(t),\sigma(t)),\quad\gamma,\sigma\in W_{o}(M).$$

In order to construct Dirichlet forms  associated to stochastic heat equations on Riemannian path space, we need to introduce the following $L^1$-distance, which is a smaller distance than the above uniform distance $d_\infty$ on  $W_{o}(M)$:
$$\tilde{d}(\gamma,\eta):={\int_0^1\rho(\gamma_s,\eta_s)\dd s}, \quad \gamma,\eta\in W_o(M).$$  Let $E$  denote the closure of $W_o(M)$
in $\{\eta:[0,1]\rightarrow M; \int_0^1\rho(o,\eta_s)ds<\infty\}$ with respect to the distance $\tilde{d}$. Then $E$ is a Polish space.

Let $O(M)$ be the orthonormal frame bundle over $M$, and let $\pi: O(M) \rightarrow M$ be the canonical projection.
Choosing a standard othornormal basis $\{H_i\}_{i=1}^d$
of horizontal vector fields on $O(M)$, and consider the following SDE,
\begin{equation}\label{eq2.1}
\begin{cases}
&d U_t=\sum^d_{i=1}H_i(U_t)\circ \dd B_t^i,\quad t\geq0\\
& U_0=u_o,
\end{cases}
\end{equation}
where $u_o$ is a fixed orthonormal basis of $T_o M$ and $B^1_t,\cdots,B_t^d$ are independent
Brownian motions on $\mathbb{R}$.
Then  $x_t:=\pi(U_t),\ t\geq0$ is the Brownian motion on
$M$ with initial point $o$, and $U_{\cdot}$ is the (stochastic) horizontal lift along
$x_{\cdot}$. Let $\mu_o$ be the distribution of $x_{[0,1]}$, then $\mu_o$
is a probability measure on $W_{o}(M)$.

Let $\F C_b^1$ be the space of bounded Lipschitz
continuous cylinder functions on $W_o(M),$ i.e. for every $F\in \F C_b^1$, there exist some $m\geq1,~g_i\in \mathrm{Lip}(M), ~m\in \mathbb{N},~ f\in C_b^1(\mathbb{R}^m)$
such that
\begin{equation}\label{eq2.2}\aligned
F(\gamma)=f\left(\int_0^1  g_1(s,\gamma_s) \dd s,\int_0^1  g_2(s,\gamma_s) \dd s,...,\int_0^1  g_m(s,\gamma_s) \dd s\right),\quad \gamma\in W_o(M),\endaligned\end{equation}
where
$$\mathrm{Lip}(M):=\left\{g:[0,1]\times M\rightarrow\mathbb{R},|g(s,\eta_s)-g(s,\gamma_s)|\leq C \rho(\eta_s,\gamma_s),s\in[0,1], \eta,\gamma\in E\right\}.
$$
For any $F\in \F C_b^1$ with \eqref{eq2.2} form and $h\in \HH:=L^2([0,1];\mathbb{R}^d)$, the directional derivature of $F$ with respect to $h$ is given by
 $$D_hF(\gamma)=\sum_{j=1}^m\hat{\partial}_jf(\gamma)\int_0^1
 \left\langle U_s^{-1}(\gamma)\nabla g_j(s,\gamma_s),h_s
 \right\rangle_{\mathbb{R}^d} \dd s,\quad \gamma\in W_o(M),$$
where
$$\hat{\partial}_jf(\gamma):=\partial_jf\bigg(\int_0^1  g_1(s,\gamma_s) \dd s,\int_0^1  g_2(s,\gamma_s) \dd s,...,\int_0^1  g_m(s,\gamma_s) \dd s\bigg),$$
and for $\gamma\in E\backslash W_o(M)$ we define $D_hF(\gamma)=0$. By Riesz's representation theorem,
there exists a gradient operator $DF(\gamma)\in \HH$ such that $\langle DF(\gamma),h\rangle_{\HH}=D_hF(\gamma), \gamma\in E, h\in \HH$. In particular, for $\gamma\in W_o(M)$,
$DF(\gamma)=\sum_{j=1}^m\hat{\partial}_jf(\gamma)U_s^{-1}(\gamma)\nabla g_j(s,\gamma_s).$ We call $DF$ the ${\bf L^2}$-gradient of $F$ on path space.
Denote by $\H$ the Cameron-Martin space:
$$\mathbb{H}:=\left\{h\in C^1([0,1];\mathbb{R}^d)\Big| h(0)=0, \|h\|^2_\H:=\int^1_0\|h'(s)\|^2\dd s<\infty\right\}.$$
 Taking $\{e_k\}\subset \H$ such that it is an orthonormal basis in $\HH$, consider the following symmetric quadratic form
$$\E(F,G):=\frac{1}{2}\int_E \langle DF, DG\rangle_{\HH}\dd\mu_o=\frac{1}{2}\sum_{k=1}^\infty\int_E D_{e_k}F D_{e_k}G  \dd\mu_o; \quad F,G\in\F C_b^1.$$

\beg{thm}\label{T2.1}  The quadratic form $(\E, \F C_b^1)$
is closable and its closure $(\E,\D(\E))$ is a quasi-regular Dirichlet form on $L^2(E;\mu_o)=L^2(W_o(M);\mu_o)$.
\end{thm}
\noindent \emph{Sketch of the proof}: For the compact Riemannian manifold, we can derive the closability of  $(\E, \F C_b^1)$ by the integration by parts formula in \cite{DR92} along each $e_k$. By a localization technique, the integration by parts formula in \cite{DR92} also can be extended to the general Riemannian manifolds, which implies the closability in the general case. The quasi-regularity of the Dirichlet form follows essentially by the same argument as in \cite{MR92}.$\hfill\Box$

\vskip.10in

By using the theory of Dirichlet forms (refer to \cite{MR92}), we obtain:

 \beg{thm}\label{T2.2} There exists a conservative (Markov) diffusion process
 $M=(\Omega,\F,\M_t,$ $(X(t))_{t\geq0},(\mathbf{P}^z)_{z\in E})$ on $E$ \emph{properly associated with} $(\E,\D(\E))$, i.e. for $u\in L^2(E;\mu_o)\cap\B_b(E)$,
 the transition semigroup $P_tu(z):=\mathbf{E}^z[u(X(t))]$ is a $\E$-quasi-continuous version of $T_tu$ for all $t >0$,
 where $T_t$ is the semigroup associated with $(\E,\D(\E))$.
\end{thm}

 Here for the notion of $\E$-quasi-continuity we refer to \cite[ChapterIII,Definition 3.2]{MR92}. By Fukushima's decomposition we have
\beg{thm}\label{T2.3}  There exists a \emph{properly  $\E$-exceptional set} $S\subset E$, i.e. $\mu_o(S)=0$ and $\mathbf{P}^z[X(t)\in E\setminus S, \forall t\geq0]=1$ for $z\in E\backslash S$, such that $\forall z\in E\backslash S$ under $\mathbf{P}^z$,  the sample paths of the associated  process $M=(\Omega,\F,\M_t,$ $(X(t))_{t\geq0},(\mathbf{P}^z)_{z\in E})$ on $E$ satisfy the following: for $u\in \D(\E)$
 \begin{equation}\label{eq2.3}\aligned u(X_t)-u(X_0)=M_t^u+N_t^u\quad \mathbf{P}^z-a.s.,\endaligned\end{equation}
 where $M^u$ is a martingale with quadratic variation process given by $\int_0^t |Du(X_s)|_\HH^2ds$ and $N^u$ is a zero quadratic variation process.
 In particular, for $u\in D(L)$, $N_t^u=\int_0^tLu(X_s)ds$, where $L$ is the generator of $(\E,\D(\E))$.
\end{thm}

\beg{Remark}\label{r2.4}

\beg{enumerate}
\item[$(a)$]  If we choose $u(\gamma)=\int_{r_1}^{r_2}u^\alpha(\gamma_s)ds\in\mathcal{F}C_b^1$, with local coordinates $u^\alpha$   on $M$, then the quadratic variation process for $M^u$ is the same as  that for the martingale part in (1.1).

\item[$(b)$]  Theorems 2.2-2.3 still hold if the path space is replaced by the loop space(or the free path and free loop cases) and Wiener measure is  replaced by the associated measure under some suitable conditions.
\end{enumerate}
\end{Remark}

\section{Properties of SHE}\label{sect3}
In this section, we will study properties of $X_t,t\geq0,$ constructed in Section 2. First we present the logarithmic Sobolev inequality for the damped gradient $\tilde{D}F$ assuming $M$ is stochastically complete, which implies
the  logarithmic Sobolev inequality for the Dirichlet form considered in Section 2.

For any $F\in \F C_b^1$, we define the damped gradient $\tilde{D}F$ of $F$ by
$$\tilde{D}F(t)=M_t^{-1}\int_t^1M_s(DF(s))\dd s,$$
where $M_t$ is the solution of the equation
$$\frac{\dd}{\dd t}M_t+\frac{1}{2}M_tRic_{U_t}=0,\quad M_0=I.$$
Suppose that $\Ric\geq -K$ for  $K\in\mathbb{R}$. Define the quadratic form corresponding to $\tilde{D}F$ by
$$\tilde{\E}(F,G)=\frac{1}{2}\int_E \langle \tilde{D}F, \tilde{D}G\rangle_{\HH}\dd\mu_o,\quad F,G\in\F C_{b}^1.$$

 \beg{thm}\label{T3.1}[Log-Sobolev inequality] Suppose that $\Ric\geq -K$ for  $K\in\mathbb{R}$. The log-Sobolev inequality holds for
$(\tilde{\E},\D(\tilde{\E}))$, i.e.,
\begin{equation*}\label{e16}
\mu_o(F^2\log F^2)\le 2 \tilde{\E}(F,F),\ \ \ \ F \in \F C^1_{b},
\ \mu_o(F^2)=1.
\end{equation*}
In particular, we have
$$\mu_o(F^2\log F^2)\le 2C(K) \E(F,F),\ \ \ \ F \in \F C^1_{b},
\ \mu_o(F^2)=1$$
where $C(K)=\frac{e^K-1-K}{K^2}\wedge C_0(K)$ with
\begin{equation*}
C_0(K)=\begin{cases}
& \frac{4}{K^2}\left(1-\sqrt{2e^{\frac{K}{2}}-e^{K}}\right),\quad\quad\text{if}~K<0,\\
&\frac{2}{K^2}\left(e^{K}-2e^{\frac{K}{2}}+1\right),\quad\quad\quad\text{if} ~K>0.
\end{cases}
\end{equation*}
\end{thm}

\beg{Remark}\label{r3.2}  \beg{enumerate}
\item[$(i)$]  In fact, Theorem 3.1 had first been proved in \cite{GW06}. Compared to the results in there, our constant $C(K)$ is smaller. By comparing the classical O-U Dirichlet form
and the $L^2$-Dirichlet form, we note that the the LSI associated to the two Dirichlet forms are essentially different, the former requires
upper and lower bounds of the Ricci curvature of $M$, and the latter only needs a lower bound for the Ricci curvature.

\item[$(ii)$] According to \cite{W05}, the log-Sobolev inequality implies hypercontractivity of the associated semigroup $P_t$, in particular, the  $L^2$-exponential ergodicity of the process:
$\|P_tf-\int f d\mu_o\|_{L^2}^2\leq e^{-t/C(K)}\|F\|_{L^2}^2.$

\item[$(iii)$] The log-Sobolev inequality  also implies the irreducibility of the Dirichlet form $(\E,\D(\E))$. It is obvious that the Dirichlet form $(\E,\D(\E))$ is recurrent.
Combining these two results, by [FOT94, Theorem 4.7.1], for every nearly Borel non-exceptional set $B$,
$$\mathbf{P}^x(\sigma_B\circ\theta_n<\infty,\forall n\geq0)=1, \quad \textrm{ for q.e. } x\in X.$$
Here $\sigma_B=\inf\{t>0:X_t\in B\}$, $\theta$ is the shift operator for the Markov process $X$, and for the definition of  nearly Borel non-exceptional set we refer to [FOT94]. Moreover by [FOT94, Theorem 4.7.3] we obtain the following strong law of large numbers: for $f\in L^1(E,\mu_o)$
$$\lim_{t\rightarrow\infty}\frac{1}{t}\int_0^tf(X_s)\dd s=
\int f\dd\mu_o, \quad \mathbf{P}^{x}-a.s.,$$
for q.e. $x\in E$.

\end{enumerate}
\end{Remark}

\noindent \emph{Sketch of the proof of Theorem 3.1}: The proof follows from the following martingale representation: for
$F\in L^2(\mu_o)$,
$$F=\EE(F)+\int^1_0\left\langle \EE\left[M_s^{-1}\int_s^1M_\tau(DF(\tau))\dd\tau\bigg|\F_s\right], \dd W_s\right\rangle ,$$
and some delicate estimates. Here $W$ is the anti-development of $\gamma$ and $\{\F_s\}$ is the filtration generated by $W$.$\hfill\Box$
\vskip.10in

Upper and lower bounds of the Ricci curvature on a Riemannian manifold were well characterized by the diffusion process associated to the O-U Dirichlet form given by the
Malliavin gradient in \cite{N}. If the O-U Dirichlet form is replaced by our $L^2$-Dirichlet form, then we can only obtain the
following  characterizations for the lower bound of the Ricci curvature. This further indicates that these two processes have essential differences.

In fact, the results in Section 2 and Theorem 3.1 also hold when we change $1$ to any $T>0$. To state our results, let us first introduce some notations:
For any point $y\in M$ and $T>0$, let $x_{y,[0,T]}$ be the Brownian motion starting from $y\in M$ up to time $T$, and $\mu_{T,y}$ be the distribution of Brownian motion $x_{y,[0,T]}$ on $W_y^T(M):=\{\gamma\in
C([0,T];M)|\gamma(0)=y\}$.
For any $n\geq1$ and $G\in \F C_b^T$ with $\F C_b^T$ defined as in (2.2) with $1$ replaced by $T$, define
$$\aligned\mathcal{E}^K_{T,n,y}(G,G)&=(1+n)C_1(K)\int_{W_y^T(M)}\int_0^{T-\frac{1}{n}}|DG(\gamma)(s)|_{\mathbb{R}^d}^2\dd s\dd  \mu_{T,y}(\gamma)\\
&~~~~~+\left(\frac{1}{n}+\frac{1}{n^2}\right)C_{2,n}(K)\int_{W^T_y(M)}\int_{T-\frac{1}{n}}^T|DG(\gamma)(s)|_{\mathbb{R}^d}^2\dd s\dd\mu_{T,y}(\gamma).\endaligned$$
where
$$C_1(K)=\left[\frac{1}{K^2}\big(TKe^{KT}-e^{KT}+1\big)\right]\bigvee \frac{T^2}{2}, \quad C_{2,n}(K)=\frac{e^{KT}-1}{K}\left(1\vee e^{-\frac{K}{n}}\right).$$
Let $p_t$ be the Markov semigroup of the process  $x_y$ given by $p_tf(y)=\EE [f(x_{t,y})],y\in M,f\in \B_b(M),t\geq0$.
Denote by $C_0^\infty(M)$ the set of all smooth  functions with compact support on $M$.
\vskip.10in

 \beg{thm}\label{T3.3}  For  $K\in\mathbb{R}$, the following statements are equivalent:
\beg{enumerate}
\item[$(1)$]
$\Ric\geq -K $.

\item[$(2)$] For every $f\in C_0^\infty(M),T>0$ and $y\in M$, we have
$$\left|\int^T_0\nabla p_sf(y)d s\right| \leq \int^T_0\e^{\frac{Ks}{2}} p_s|\nabla f|(y) \dd s .$$

\item[$(3)$]
 For every $y\in M,T>0$, the following log-Sobolev inequality holds for every $n\in\mathbb{N}$:
$$\mu_{T,y}(F^2\log F^2)\leq 2\E^K_{T,n,y}(F,F),\quad F\in\F C^T_b,~\mu_{T,y}(F^2)=1.$$

\item[$(4)$] For every $y\in M,T>0$, the following Poincar\'{e}-inequality holds for every $n\in\mathbb{N}$:
$$\mu_{T,y}(F^2)\leq \E^K_{T,n,y}(F,F),\quad F\in\F C^T_b,~\mu_{T,y}(F)=0.$$

\end{enumerate}
\end{thm}

\noindent \emph{Sketch of the proof}: 1)$\Longrightarrow$ 2) follows from the gradient formula. Conversely, taking $F(\gamma):=\int^T_0f(\gamma_s)ds$ for some function $f\in C^1_0(M)$ with
\begin{equation}\label{eq3.1}f\in C_0^\infty(M),\quad|\nabla f|(y)=1,\quad \textrm{Hess}_f(y)=0,\end{equation} and applying $F$ into $2)$,  $1)$ can be derived from  the following formula in \cite{W14}$$\frac{1}{2}\Ric(\nabla f,\nabla f)(y)=\lim_{T\downarrow0}\frac{p_T|\nabla f|(y)-|\nabla p_Tf|(y)}{T}.$$

1)$\Rightarrow$ 3) follows similarly as in the proof of Theorem \ref{T3.1}.

3)$\Rightarrow$ 4) is standard.

4)$\Rightarrow$ 1): For each $k\geq1$, take $F(\gamma)=k\int_{T-1/k}^Tf(\gamma_s)ds$ for some $f$ as \eqref{eq3.1}. Then using  this formula
$$\frac{1}{2}\Ric(\nabla f,\nabla f)(y)=\lim_{T\downarrow0}\frac{1}{T}\bigg(\frac{p_Tf^2(y)-(p_Tf)^2(y)}{2T}-|\nabla p_Tf(y)|^2\bigg),$$ it is not difficult to obtain $1)$.
$\hfill\Box$

\beg{thebibliography}{99}

\leftskip=-2mm
\parskip=-1mm

\bibitem{AD99}  L. Andersson,  B. K. Driver, \emph{Finite-dimensional approximations to
Wiener measure and path integral formulas on manifolds}, J. Funct. Anal, 165, no. 2, (1999), 430–498.

\bibitem{ALR93} S. Albeverio, R. L\'{e}andre and M. R\"{o}ckner: Construction of a rotational invariant diffusion on
the free loop spaces, C.R. Acad. Sci. Paris, 316, Serie I (1993), 287-292.

\bibitem{BHZ16}Y. Bruned, M. Hairer, and L. Zambotti. \emph{Algebraic renormalisation of regularity structures.
arXiv:1610.08468,} pages 1-84, 2016.

\bibitem{CH16}A. Chandra and M. Hairer. \emph{An analytic BPHZ theorem for regularity structures}
arXiv:1612.08138, pages 1-113, 2016.

\bibitem{CW14}X. Chen, B. Wu, Functional inequality on path space over a non-compact Riemannian
manifold, J. Funct. Anal. 266(2014), 6753-6779.

\bibitem{Dri92}B. K. Driver, \emph{A Cameron-Martin quasi-invariance theorem for Brownian motion on
a compact Riemannian manifolds,} J. Funct. Anal. 110 (1992) 273–376.

\bibitem{DR92}B. K. Driver and M. R\"{o}ckner,\emph{Construction of diffusions on path and loop spaces of
compact Riemannian manifolds,} C. R. Acad. Sci. Paris Ser. I 315 (1992) 603–608

\bibitem{Fun92} T. Funaki, \emph{A stochastic partial differential equation with values in a manifold,}
J. Funct. Anal. 109, no. 2, (1992), 257-288.

\bibitem{FOT94}M. Fukushima, Y. Oshima, and M. Takeda, \emph{Dirichlet Forms and Symmetric
Markov Processes,} de Gruyter, Berlin (1994), second edition (2010)

\bibitem{GW06} M. Gourcy,  L. Wu, \emph{Logarithmic Sobolev inequalities of diffusions for the L2-metric}. Potential
 Anal. 25 77-102 (2006)

\bibitem{Hai14} M. Hairer, \emph{A theory of regularity structures,} Invent. Math., 198(2):269-504, 2014

\bibitem{Hai16}M. Hairer,\emph{The motion of a random string,} arXiv:1605.02192, pages 1-20, 2016.

\bibitem{MR92}Z. M. Ma and M. R\"{o}ckner, \emph{Introduction to the Theory of (Non-Symmetric) Dirichlet
Forms}, Springer-Verlag, Berlin, Heidelberg, New York, 1992.

\bibitem{N} A. Naber, Characterizations of bounded Ricci curvature on smooth and nonsmooth spaces, {\it arXiv: 1306.6512v4}.

\bibitem{Nor98} J. R. Norris,\emph{Ornstein-Uhlenbeck processes indexed by the circle}, Ann.
Probab. 26, no. 2, (1998), 465-478.

\bibitem{RWZZ17}M. R\"{o}ckner, B. Wu, R. Zhu, X. Zhu,  in preparation

\bibitem{W05} F.Y. Wang, \emph{Functional Inequalities, Markov Semigroup and Spectral Theory}. Chinese Sciences
Press, Beijing (2005)

\bibitem{W14}F.- Y. Wang, Analysis for diffusion processes on Riemannian manifolds, World
Scientific, (2014)

\bibitem{WW1} F.- Y. Wang and B. Wu, \emph{Quasi-Regular Dirichlet Forms on
Riemannian Path and Loop Spaces,} Forum Math. 20(2008), 1085--1096.

\end{thebibliography}
\end{document}